\documentclass[11pt]{amsart}

\usepackage{graphicx}
\usepackage{color}
\usepackage{amssymb}
\usepackage{latexsym}

\newcommand{\R}{\mathbb{R}}

\newcommand{\be}{\begin{equation}}
\newcommand{\ee}{\end{equation}}

\newtheorem{thm}{Theorem}[section]
\newtheorem{lem}[thm]{Lemma}

\newtheorem{pro}[thm]{Proposition}

\author[O. Goubet]{O. Goubet}
\author[L. Molinet]{L. Molinet}

\title[Nonlinear Schr\"odinger Equations]{Global attractor for weakly damped Nonlinear
Schr\"odinger equations in $L^2(\R)$}
\address[Olivier Goubet]{LAMFA CNRS UMR 6140 \\
Universit\'e de Picardie Jules Verne \\
33 rue Saint-Leu 80039 Amiens cedex.}
\address[Luc Molinet] {LAGA, Institut Galil\'ee, Universit\'e Paris 13 \\
93430 Villetaneuse}

\email{olivier.goubet@u-picardie.fr}

\begin{document}

\begin{abstract}
We prove that the weakly damped nonlinear Schr\"odinger flow in
$L^2(\mathbb{R})$ provides a dynamical system which possesses a
 global attractor. The proof relies on the continuity of
the Schr\"odinger flow for the weak topology in $L^2(\R)$.
\end{abstract}


\subjclass[2000]{}

\date{23 march 2007}

\maketitle

\section{Introduction}

Nonlinear Schr\"odinger (NLS) and Korteweg-de Vries equations are
asymptotical models for the waterwave propagation. These models
supplemented with a damping and an external force provide examples
of infinite-dimensional dynamical systems, in the framework
described in \cite{Temam}, \cite{Hale}, \cite{Raugel}, \cite{MZ}. We
focus here on the cubic NLS equation

\begin{equation}\label{equation}
u_t+\gamma u +iu_{xx}+i|u|^2u=f,
\end{equation}
\noindent where $\gamma>0$ is the damping parameter and where the
external forcing $f(x)$, that is independent of $t$, belongs to
$L^2(\R)$.

To define an infinite-dimensional dynamical system from this
evolution equation, we supplement (\ref{equation}) with an initial
data $u_0$ in some Sobolev space $X$ (or in some complete metric
space) such that the corresponding initial value problem is well
posed (in the Hadamard sense: existence and uniqueness of
trajectories $u(t)=S(t)u_0$ in $X$, continuity of $S(t):u_0\mapsto
u(t)$ in $X$). This short article is concerned with the existence
of a global attractor of the NLS flow with low regular initial
data in $L^2(\R)$. Recall that a global attractor is a compact
set, invariant by the flow,  that attracts all trajectories
uniformly on bounded sets. Note that the existence of a such set
for the NLS equations in  more regular function spaces, as for
instance $H^1(\R)$, is
 well known (see \cite{Akroune} and the references
therein).

Our main result states as follows

\begin{thm}\label{attractor}
The semi-group $S(t)$ provides an infinite-dimensional dynamical
system in $L^2(\R)$ that has a global attractor $\mathcal{A}$.
\end{thm}

Let us describe the strategy of the proof. As in \cite{Gh1}, we
first prove that the NLS flow features a weak global attractor in
$L^2(\R)$, that is a global attractor for the weak topology of
$L^2(\R)$. For that purpose, we need to establish the following
result that is interesting on its own :

\begin{thm}\label{weakcontinuity}
The semi-group $S(t)$ is a continuous mapping for the weak topology
of $L^2(\R)$.
\end{thm}

\noindent The usual arguments to prove this kind of result in
 a function space $ X $ make use of the fact that the initial value problem is
 well-posed in a function space where $ X $ is locally compactly
 embedded (cf. \cite{Gh1}). Here such an argument can not be
 invoked since  there is no available result on the
 well-posedness of the (\ref{equation}) in a space where $ L^2(\R)
 $ is locally compactly  embedded (cf. \cite{KPV}).
Note that the situation  for KdV equations is easier according to
this last point (see \cite{GR}, \cite{Tsugawa}).

  Therefore
this result  is new and can be outlined as follows. Calling
$U(t)u_0$
 the solution of
 \be\label{linear} u_t+iu_{xx}=0; \; u(0)=u_0,
\end{equation}
 it is well-known that the {\it linear}
Schr\"odinger equation features the so-called Kato's smoothing
effect that reads:
\be\label{kato}
||D_x^{1/2}U(t)u_0||_{L^\infty_xL^2_t}\leq c||u_0||_{L^2_x}.\ee

\noindent Above,  $D_x=\sqrt{-\Delta}$ stand for the operator
with Fourier  symbol $|\xi|$.
 Using the Christ-Kiselev theorem \cite{CK} (as in
\cite{MolinetRibaud} in another context), we are able to prove
that this smoothing effect is also valid for the nonlinear
Schr\"odinger equation. Then the weak continuity is valid due to
some compactness argument that allow us to pass to the limit in
the nonlinear term.

\noindent These arguments are developed in Section 1 below. In
Section 2 we complete the proof of Theorem \ref{attractor}. First
we prove the existence of a weak attractor. Then, using the
famous J. Ball's argument (see \cite{Ball}, \cite{Wang},
\cite{MRW}), we establish that the weak attractor is actually a
global attractor in the usual sense.

\section{Continuity of the flow for the weak topology}

To begin with, we observe that for finite time results the damping
parameter and the external forcing do not play a role. Then in this
section we may assume for the sake of simplicity that $\gamma=0$ and
$f=0$.

The usual way to solve the IVP problem associated to

\be\label{1} u_t+iu_{xx}+i|u|^2u=0,\ee

\noindent supplemented with initial data $u_0$ is to perform a fixed
point argument for the Duhamel's form of (\ref{1}) that reads

\be\label{2} u(t)=U(t)u_0-i\int_0^t U(t-s)|u(s)|^2u(s)ds. \ee

\noindent Thanks to well-known Strichartz inequalities, we usually
perform a fixed point into the space $C([0,T],L^2(\R))\cap L^6_{T,x}$,
where $L^6_{T,x}=L^6([0,T]\times\R_x)$. We first state and prove

\begin{pro}
There exists a numerical constant $c$ such that for a solution to
(\ref{2}) \be\label{3}||D_x^{1/2}u||_{L^\infty_xL^2_T}\leq
c(||u_0||_{L^2_x}+||u_0||_{L^2_x}^3).\ee
\end{pro}

\noindent Proof of the Proposition: the key point is to estimate
the nonlinear term in (\ref{2}). For that purpose, we first recall
the dual estimate to (\ref{kato}) that reads

\be\label{4}||\int_\R U(-s)D_x^{1/2}Gds||_{L^2_x}\leq
c||G||_{L^1_xL^2_t}. \ee

\noindent We now prove

\be\label{5}||\int_{\R_s} U(t-s)D_x^{1/2}fds||_{L^\infty_xL^2_t}
\leq c||f||_{L^\frac65_{t,x}}. \ee \noindent  Actually, following
P. Tomas duality argument, it is equivalent to prove that for any
smooth function $G$ that satisfies $||G||_{L^1_xL^2_t}\leq 1$, it
holds

\be\label{6} \Bigl|\int_{\R^3}
U(t-s)D_x^{1/2}f(s,x)\overline{G}(t,x)dtdxds \Big| \leq
c||f||_{L^{\frac{6}{5}}_{t,x}}. \ee Note that  the left-hand side
member of the above estimate can be rewritten as \be\label{6}
\Bigl|\int_{\R} \Bigl( \int_{\R} U(-s)D_x^{1/2}f(s,x)\, ds \Bigr)
\Bigl(\int_{\R}\overline{U(-t)G(t,x)}dt\Bigr) \, dx  \Big| \; ,
\ee
 \noindent  Hence, applying Cauchy-Schwarz in $  x$ and using (\ref{4}),
  it finally suffices to check that

\be\label{7}||\int_\R U(-s)fds||_{L^2_x}\leq
c||f||_{L^{\frac{6}{5}}_{t,x}}, \ee

\noindent Since this is nothing else but the dual form of  the
classical Strichartz estimate for the Schr\"odinger group on $ \R
$:
 \be
 \|U(t) u_0\|_{L^6_{t,x}} \le c \|u_0\|_{L^2_x}\; ,\ee
 we are done.

\noindent Recall now from \cite{MolinetRibaud}

\begin{lem}(Christ-Kiselev)
Consider a linear operator defined on space-time functions $f(t,x)$
by \be\label{8} Tf(t)=\int_{\R_s}K(t,s)f(s)ds.\ee Assume
\be\label{9}||Tf||_{L^\infty_xL^2_t}\leq
c||f||_{L^{\frac{6}{5}}_{t,x}}, \ee \noindent then

\be\label{91}||\int_0^t K(t,s)f(s)ds||_{L^\infty_xL^2_t}\leq
c||f||_{L^{\frac{6}{5}}_{t,x}}. \ee

\end{lem}

\noindent According to \cite{MolinetRibaud}, this is valid since
$\min(+\infty,2)> {\rm max}(\frac65,\frac65)$.

\noindent We then apply this argument to the nonlinear term in
(\ref{2}), for $t\in[0,T]$. This leads to

\be\label{10}\begin{split}||\int_0^t
D_x^{1/2}U(t-s)\bar{u}^2uds||_{L^\infty_xL^2_T}\leq
c||u^3||_{L^{\frac{6}{5}}_{T,x}}\\
\leq c||u||^3_{L^{\frac{18}{5}}_{T,x}} \leq c ||u||_{L^{2}_{T,x}}
||u||^2_{L^{6}_{T,x}}\; .\end{split}\ee

\noindent We conclude the proof of the proposition using that $u$ is
bounded in $C([0,T],L^2(\R))\cap L^6_{T,x}$. $\square$

At this stage we complete the proof of Theorem
\ref{weakcontinuity}. Consider $u_{0,\varepsilon} \rightharpoonup
u_0$ in $L^2_x$. Due to the previous proposition, we know that,
for any $K$ compact subset of $\R_x$, the sequence
$u_{\varepsilon}$ remains in a bounded set of
$C([0,T],L^2(\R))\cap L^6_{T,x}\cap L^2_TH^\frac12_x(K)$. Going
back to the equation, we observe that $\partial_t
u_{\varepsilon}$ remains in a bounded set of $L^2_TH^{-2}_x$.
Hence, due to a standard compactness argument, the sequence
$u_{\varepsilon}$, up to a subsequence extraction, converges
towards some function $ v $ strongly in $L^2_TL^{2}(K)$. By
interpolation, the strong convergence is also valid in
$L^4_TL^{4}(K)$. This allows us to pass to the limit in the
equation and to conclude that the limit $ v $ is a solution of
(\ref{1})  belonging
 to the class of uniqueness
  $ L^6_{T,x} $. Set $(.,.)$ for the $L^2_x$ scalar product. By (\ref{1}) and 
  the bounds above, it is easy to check that, for any smooth
 space function $\phi $ with compact support, the family $
 \{t\mapsto (u_\varepsilon(t), \phi)\} $ is uniformly
 equi-continuous on $ [0,T] $. Ascoli's theorem then ensures that
 $ (u_\varepsilon(\cdot),\phi)$ converges to $ (v(\cdot),\phi) $
 uniformly  on $ [0,T] $ and thus $v(0)=u_0$. By uniqueness, it
 follows that $ v\equiv u $ and from the above convergence result,
   it results that $u_\varepsilon(t) \rightharpoonup u(t) $ in $L^2_x$ for all
  $ t\in [0,T] $.

$\square$

\section{Proof of the main Theorem }

To begin with, we prove the existence of an absorbing ball for the
semi-group; multiplying (\ref{equation}) by $\bar{u}$ and
integrating in $x$ the real part of the resulting equation

\begin{equation}\label{absorbing}
\frac{1}{2}\frac{d}{dt} ||u||^2_{L^2_x}+\gamma ||u||^2_{L^2_x} =
{\rm Re} \; \int f\bar{u}dx\leq \frac{1}{2}\gamma ||u||^2_{L^2_x}+
\frac{1}{2\gamma}||f||^2_{L^2_x}.
\end{equation}

\noindent This implies

\be\label{21}||u(t)||^2_{L^2_x}\leq e^{-\gamma t}||u_0||^2_{L^2_x}+
\frac{1-e^{-\gamma t}}{\gamma^2}||f||^2_{L^2_x}. \ee

\begin{pro}
 The ball $X$ of radius $M_0=2\frac{||f||_{L^2_x}}{\gamma}$ is an
absorbing set for the dynamical system under consideration.
\end{pro}

\noindent We endow then this absorbing ball with the weak topology
of $L^2_x$. $X$ is then a compact metric space and $S(t)$ acts
continuously on $X$ according to Theorem \ref{weakcontinuity}.
Therefore, using Theorem I.1.1 in \cite{Temam} the $\omega$-limit
set $\mathcal{A}=\cup_{s>0}\cap_{t>s} \overline{S(t)X}$ is a global
attractor. In fact

\be\label{22}\mathcal{A}=\{a\in X; \exists b_n \in X, t_n\rightarrow
+\infty, \; S(t_n)b_n\rightharpoonup a.\} \ee

\noindent We plan to transform this weak convergence into a strong
convergence. We use the famous J. Ball's argument. We begin with the
energy equation that asserts that for any $\tau>0$, due to
(\ref{absorbing}),

\be\label{23} ||S(t_n)b_n||^2_{L^2_x}=e^{-2\gamma
\tau}||S(t_n-\tau)b_n||^2_{L^2_x} -2{\rm Re}\int_0^\tau
\int_{\R_x}e^{-2\gamma s}\overline{f}(x)S(t_n-s)b_ndsdx. \ee

\noindent According to the weak convergence, we have

\be\label{24}\lim_{n\rightarrow +\infty}2{\rm Re}\int_0^\tau
\int_{\R_x}e^{-2\gamma s}\overline{f}(x)S(t_n-s)b_ndsdx=2{\rm
Re}\int_0^\tau \int_{\R_x}e^{-2\gamma s}\overline{f}(x)S(-s)adsdx
\; . \ee

\noindent Using once again the energy equality (\ref{absorbing}) we
also have that

\be\label{25} ||a||^2_{L^2_x}=e^{-2\gamma
\tau}||S(-\tau)a||^2_{L^2_x} -2{\rm Re}\int_0^\tau
\int_{\R_x}e^{-2\gamma s}\overline{f}(x)S(-s)adsdx. \ee

\noindent Therefore

\be\label{26}\limsup_n||S(t_n)b_n||^2_{L^2_x}\leq ||a||^2_{L^2_x}
+2e^{-2\gamma \tau}M_0^2,\ee

\noindent since for $n>\tau$  $S(t_n-\tau)b_n$ is in $X$ and
$S(-\tau)a$, that belongs to the weak attractor, remains trapped in
$X$. Letting $\tau\rightarrow +\infty$ implies that $\mathcal{A}$
attracts the bounded sets for the $L^2_x$ strong topology. To prove
that $\mathcal{A}$ is compact is very similar and then omitted.
$\square$


\begin{thebibliography}{99}

\bibitem{Akroune} N. Akroune, Regularity of the attractor for a weakly damped
Schr\"odinger equation on $I\!\!R$, {\it Applied Math. Letters}, 12,
(1999), 45-48

\bibitem{Ball} J. Ball, {\sl Global attractors for damped semilinear
 wave equations. Partial differential equations and applications},
  Discrete Contin. Dyn. Syst.  10  (2004),  no. 1-2, 31--52.

\bibitem{CK} M. Christ and A. Kiselev, {\sl Maximal functions
associated to filtrations}, J. Funct. Analysis 179, (2001), 406-425.

\bibitem{Gh1} J.-M. Ghidaglia, {\sl Finite dimensional behavior for
the weakly damped driven \par \noindent Schr\"odinger equations},
Ann. Inst. Henri Poincar\'e, 5, (1988), 365-405.

\bibitem{G1} O. Goubet, {\sl Regularity of the attractor for the weakly
damped nonlinear \par\noindent  Schr\"odinger equations}, Applicable
Anal., 60, (1996), 99-119. \par

\bibitem{GR} O. Goubet and R. Rosa, {\sl Asymptotic smoothing and
the global attractor for a weakly damped kdv equation on the real
line},  J. of Diff. Eq., vol 185, 1, pp 25-53, 2002.

\bibitem{Hale} J. Hale, {\sl Asymptotic behavior of Dissipative Systems},
Math. surveys and Monographs, vol 25, AMS, Providence, 1988.

\bibitem{KPV}  C. Kenig, G. Ponce and L. Vega, {\sl On the ill-posedness of some
canonical dispersive equations},  Duke Math. J.  106  (2001),  no.
3, 617--633.

\bibitem{MZ} A. Miranville and S. Zelik,
{\emph Attractors for dissipative partial differential equations in
bounded and unbounded domains}, Handbook of Differential Equations,
Evolutionary Partial Differential Equations, C.M. Dafermos and M.
Pokorny eds., Elsevier, Amsterdam, to appear

\bibitem{MRW} I. Moise, R. Rosa and X. Wang, {\sl Attractors for non-compact semigroups
via energy equations},  Nonlinearity  {\bf 11},  (1998),  no. 5,
1369--1393.

\bibitem{MolinetRibaud} L. Molinet and F. Ribaud, {\sl Well-posedness
results for the generalized Benjamin-Ono equation with small initial
data}, Journal de Math\'ematiques Pures et Appliqu\'ees, 83, 2004,
pp 273-311.

\bibitem{Raugel}  G. Raugel, {\sl Global attractors in partial differential equations}.
Handbook of dynamical systems, Vol. 2,  885--982, North-Holland,
Amsterdam, 2002.

\bibitem{Temam} R. Temam, {\sl Infinite Dimensional Dynamical
Systems in Mechanics and Physics,\/} Springer-Verlag, Second
Edition, 1997. \par

\bibitem{Tsugawa} K. Tsugawa, {\sl Existence of the global attractor
 for weakly damped, forced KdV equation on Sobolev spaces of negative index},
  Commun. Pure Appl. Anal.  3  (2004),  no. 2, 301--318.

\bibitem{Wang} X. Wang, {\sl An energy equation for the weakly damped
driven nonlinear Schr\"odinger equations and its applications to
their attractors}, Physica D, 88, (1995), 167-175. \par




\end{thebibliography}
\end{document}